\documentclass[english]{article}
	\usepackage[latin9]{inputenc}
	\usepackage{geometry}
	\usepackage{graphicx}
\usepackage[export]{adjustbox}
	\usepackage{amsmath}
	\usepackage{amsthm}
	\usepackage{enumerate}
	\def\disp{\displaystyle}
	\newtheorem*{theorem}{Theorem}
	\def\intro{\textbf{Introduction. }}
	\def\defi{\textbf{Definition. }}
	\def\l{\left}
	\def\r{\right}
\begin{document}
\title{An elementary proof of the rationality of $\zeta(2n)/\pi^{2n}$}
\author{Tom Moshaiov}
\date{}
\maketitle
\intro In $1735$ Euler \cite{1} proved that for each positive integer $k$, the series $\zeta(2k) = \sum_{\ell=1}^{\infty} \ell^{-2k}$ 
converges to a rational multiple of $\pi^{2k}$. 
Many demonstrations of this fact are now known, and Euler's discovery is traditionally proven using non-elementary techniques, such as Fourier series or the calculus of residues \cite{2}.  We give an elementary proof, similar to Cauchy's \cite{3} proof of the identity $\zeta(2) = \pi^2/6$, 
only extended recursively for all values $\zeta(2k)$. Our main formula 
$$\zeta(2k)=-\dfrac{(-\pi^2)^{k}}{4^{2k}-4^{k}}\left[\dfrac{4^{k}k}{(2k)!}+{\displaystyle \sum_{\ell=1}^{k-1}}(4^{2\ell}-4^{\ell})\dfrac{4^{k-\ell}}{(2k-2\ell)!}\dfrac{\zeta(2\ell)}{(-\pi^2)^{\ell}}\right] \phantom{spa}k = 1,2,3,\dots$$ 
may be derived from previously known formulae \cite{4}. 
Remarkably, Apostol \cite{5} discovered a proof similar to ours, yet arrived at a different formula, relating $\zeta(2k)$ to the Bernoulli numbers, \`a la Euler.
\\\\

\begin{theorem} Define $(a_k)_{k=1}^{\infty}$ via 
$a_{k}=(-1)^{k-1}\left(\dfrac{4^{k}k}{(2k)!}+{\displaystyle \sum_{\ell=1}^{k-1}\dfrac{(-1)^{\ell}4^{k-\ell}}{(2k-2\ell)!}a_{\ell}}\right)$. 
\footnote{We interpret the empty sum as $0$, and in paritcular $a_1=2$.} Then $\forall\phantom{.}k\ge 1$ we have
\begin{enumerate}[(i).]
\item $P_k(n)=\disp\sum_{\ell=1}^n$$\left(\cot(\frac{2\ell-1}{4n}\pi)\right)^{2k}$ equals a polynomial in $n$ of degree $2k$ with leading coefficient $a_k$.
\item $Q_k(n)=\disp\sum_{\ell=1}^n$$\left(\sin(\frac{2\ell-1}{4n}\pi)\right)^{-2k}$ equals a polynomial in $n$ of degree $2k$ with leading coefficient $a_k$.
\item $\zeta(2k)=\dfrac{\pi^{2k}}{4^{2k}-4^{k}}a_k$.\\\\
\end{enumerate}
\end{theorem}

Our proof's main ingredient is the Newton Girard formula, which relates power sums with the elementary symmetric polynomials. An elegant proof of this formula was given by Euler in \cite{6}.
\\\\

\defi Given a finite sequence $t_1,\dots,t_n$, its $r$-th {\it{power sum}} $t_1^r+\dots+t_n^r$ is denoted by $f_r$, and its $r$-th {\it{elementary symmetric polynomial}} $\displaystyle\sum_{I\in{[n]\choose r}}\prod_{i\in I}t_i$ is denoted by $e_r$, where ${[n]\choose r}$ is the family of subsets of $\{1,\dots,n\}$ of cardinality $r$. 
Note that $e_r$ is the coefficient of $t^{n-r}$ in the polynomial $(t+t_1)\dots(t+t_n)$.
\\\\

\textbf{The Newton-Girard formula. } For $k\ge 1$ we have $$f_k=(-1)^{k-1}\left[ke_k + \sum_{\ell=1}^{k-1}(-1)^\ell f_\ell e_{k-\ell}\right]$$

\textbf{Proof of theorem. }
\proof[Proof of (i).] The roots of $f(z)=\disp\sum_{k=0}^n {2n\choose 2k}z^{2k}=\frac{1}{2}\left[(z+1)^{2n}+(z-1)^{2n}\right]$ are easily verified to be 
$i\cot(\frac{2\ell-1}{4n}\pi)$, $\ell=1,\dots,2n$, 
implying $f(z)=\disp\prod_{\ell=1}^n$$(z^2+\cot^2(\frac{2\ell-1}{4n}\pi))$. Thus 
$$\tilde{P}_k(n)=e_k\l(\cot^2\l(\frac{1}{4n}\pi\r),\cot^2\l(\frac{3}{4n}\pi\r),\dots, \cot^2\l(\frac{2n-1}{4n}\pi\r)\r)={2n \choose 2k}$$ 
eqauls a polynomial in $n$ of degree $2k$ and leading coefficient $\dfrac{4^k}{(2k)!}$. 
Using induction and applying the Newton-Girard formula, we see that 
$P_k(n)=f_k(\cot^2(\frac{1}{4n}\pi),\cot^2(\frac{3}{4n}\pi),\dots, \cot^2(\frac{2n-1}{4n}\pi))$ 
equals a polynomial in $n$ of degree $2k$ with leading coefficient $a_k$. \qed
\\

\proof[Proof of (ii).] Similarly to the above, we let $\tilde{Q}_k(n)=e_k(\sin^{-2}(\frac{1}{4n}\pi),\sin^{-2}(\frac{3}{4n}\pi),\dots, \sin^{-2}(\frac{2n-1}{4n}\pi))$. 
Since $\sin^{-2}(x)=1+\cot^2(x)$, the easily verified formula $e_k(t_1+1,\dots,t_n+1)=\disp\sum_{\ell=0}^k{n-\ell\choose k-\ell}e_\ell(t_1,\dots,t_n)$ 
implies that $\tilde{Q}_k(n)$ is also a polynomial in $n$ of degree $2k$ and leading coefficient $\dfrac{4^k}{(2k)!}$. 
Again, Newton-Girard implies that $Q_k(n)=f_k(\sin^{-2}(\frac{1}{4n}\pi),\sin^{-2}(\frac{3}{4n}\pi),\dots, \sin^{-2}(\frac{2n-1}{4n}\pi))$ equals a polynomial in $n$ of degree $2k$ and leading coefficient $a_k$. \qed
\\

\proof[Proof of (iii).] For $0<x<\frac{\pi}{2}$ we have $0<\sin x<x<\tan x$. 
Therefore 
$$ P_k(n) < \disp \l(\frac{4n}{\pi}\r)^{2k} +  \l(\frac{4n}{3\pi}\r)^{2k} + \dots +  \l(\frac{4n}{(2n-1)\pi}\r)^{2k} < Q_k(n)  $$
As $P_k(n)\sim Q_k(n)\sim a_k n^{2k}$, we deduce that 
$\disp\sum_{\ell\ge 1}(2\ell-1)^{-2k}=\frac{\pi^{2k}}{4^{2k}}a_k$. 
We finish by noting that $\disp\sum_{\ell\ge 1}(2\ell-1)^{-2k}=(1-2^{-2k})\zeta(2k)$. \qed 
\\\\\\\\


\begin{thebibliography}{1}
\bibitem{1}
L. Euler, De summis serierum reciprocarum, Comment. Acad. Sci. Petropolit., 7 (1734/35), (1740) 123-134; Opera
omnia, Ser. 1 Bd. 14, 73-86. Leipzig-Berlin, 1924.
\bibitem{2}
Robin Chapman, Evaluating $\zeta(2)$. http://empslocal.ex.ac.uk/people/staff/rjchapma/etc/zeta2.pdf
\bibitem{3}
Cauchy, A. L., Bradley, R. E., \& Sandifer, C. E. (2009). Cauchy's Cours d'analyse: An annotated translation.
Dorcrecht: Springer. Note VIII.
\bibitem{4}
Saalsch\"utz, Louis. Vorlesungen \"uber die Bernoullischen Zahlen: ihren Zusammenhang mit den Secanten-
Coefficienten und ihre wichtigeren Anwendungen. Springer, 1893. pp. 14.
\bibitem{5}
Apostol, Tom M. "Another elementary proof of Euler's formula for $\zeta(2n)$." The American Mathematical Monthly
80.4 (1973): 425-431.
\bibitem{6}
Euler, Leonhard. "A double demonstration of a theorem of Newton, which gives a relation between the coefficient
of an algebraic equation and the sums of the powers of its roots." arXiv preprint arXiv:0707.0699 (2007).
\bibitem{7}
D. Zeilberger, A combinatorial proof of newton's identity, Discrete Mathematics 49 (1984), 319.
\end{thebibliography}
\end{document}